\documentclass[10pt]{article}
\usepackage{amsfonts}
\usepackage{latexsym}
\usepackage{cite}
\usepackage{amsmath,amsfonts,latexsym,amssymb}
\usepackage[mathscr]{eucal}
\usepackage{cases}
\usepackage{amsthm,color}

\usepackage[bf,small]{caption2}
\usepackage{float}
\usepackage{graphicx}
\usepackage{amsmath}
\usepackage{amssymb}
\usepackage[all]{xy}

\newtheorem{theorem}{theorem}[section]
\newtheorem{thm}[theorem]{Theorem}
\newtheorem{lem}[theorem]{Lemma}
\newtheorem{prob}[theorem]{Problem}

\newtheorem{cor}[theorem]{Corollary}

\newtheorem{nota}[theorem]{Notation}
\newtheorem{defn}[theorem]{Definition}
\newtheorem{exmp}[theorem]{Example}

\newtheorem{rmk}[theorem]{Remark}

\begin{document}

\title{\vspace{-2cm}\textbf{A new approach to the genus spectra of abelian $p$-groups}}

\author{\Large Haimiao Chen \hspace{10mm} Yang Li}

\date{}
\maketitle

\begin{abstract}
  Given a finite group $G$, the {\it genus spetrum} ${\rm sp}(G)$ of $G$ is the set of integers $g\geq 0$ such that $G$ can act faithfully on an orientable closed surface of genus $g$ by orientation-preserving homeomorphisms. The determination of ${\rm sp}(G)$ is a classical topic and has a long history, but progress is lacked. In this paper, when $G$ is an abelian $p$-group with $p>2$, we propose a new approach to ${\rm sp}(G)$, giving a structural description for ${\rm sp}(G)$ in terms of a function which can be computed in finitely many steps.

\end{abstract}

\section{Introduction}

Let $\Sigma_{g}$ denote an orientable closed surface of genus $g$.
A celebrated result of Hurwitz \cite{Hu1893} states that, if $g\ge 2$ and $G$ is a finite group acting faithfully on $\Sigma_g$ by orientation-preserving homeomorphisms, then $|G|\leq 84(g-1)$.
Wiman \cite{Wi1895} and Harvey \cite{Ha66} proved that $|G|\leq 2(2g+1)$ when $G$ is cyclic. Harvey \cite{Ha66} and Kulkarni \cite{Ku87} gave
a sharper bound for $p$-groups. These bounds motivated many works of Burnside and Maschke and others on the actions of specific kinds of groups on surfaces, and on embeddings of graphical representations of groups on surfaces (see \cite{Bu1897,Ma1896}).

Given a finite group $G$, the set ${\rm sp}(G)$ of $g\ge 0$ such that $G$ acts faithfully on $\Sigma_{g}$ by orientation-preserving homeomorphisms is called the {\it genus spectrum} of $G$. In \cite{MM92}, the determination of ${\rm sp}(G)$ was called the
{\it Hurwitz Problem} for $G$.

The genus spectrum has been studied for various classes of groups, such as cyclic $p$-groups \cite{KM91}, $p$-groups of cyclic $p$-deficiency $\leq 2$ \cite{MT98}, $p$-groups of exponent $p$ and $p$-groups of maximal class \cite{Sa09}, split metacyclic groups of order $pq$ with $p,q$ being different prime numbers \cite{We01,OW11}. The study for finite $p$-groups of maximal class in \cite{Sa09} was incomplete; it was recently finished in \cite{Sa23} and opened up a direction towards ``Co-class wise uniformity conjecture''.
Genus spectra of general abelian $p$-groups were studied in \cite{Ta11,MS19}, but are still largely unknown.


In this paper, we study ${\rm sp}(G)$ when $G$ is an abelian $p$-group with $p>2$.
Following \cite{MS19}, let $\Delta(G)=|G|/\exp(G)$, and let
$${\rm sp}_0(G)=\Big\{\frac{g-1}{\Delta(G)}\colon g\in{\rm sp}(G)\Big\}\subseteq\{-1\}\cup\mathbb{N}.$$
Call ${\rm sp}_0(G)$ the {\it reduced genus spectrum}, and denote its minimal element by $\mu_0=\mu_0(G)$. There is a {\it minimal reduced stable genus} $\sigma_0=\sigma_0(G)$ such that $g'\in{\rm sp}_0(G)$ for all $g'\ge\sigma_0$. Call the set of integers in $[\mu_0,\sigma_0]\setminus{\rm sp}_0(G)$ the (reduced) {\it spectral gap} or {\it gap sequence}. The genus spectrum is completely determined by $\mu_0,\sigma_0$ and the spectral gap. See \cite{Ku87,KM91} for similar notions in the general setting of finite groups.

The approach we take relies on a fine use of the structure of an abelian group of prime power order. We break up the Galois extension of $G$ acting on $\Sigma_g$ at the subgroup $N$ generated by ``elliptic elements", and show a nice connection of the break up with the invariants of $N$.
The main novel contribution of the paper is a clear structural description of ${\rm sp}_0(G)$ in terms of a function $\lambda$ depending on the invariants of $G$; see Theorem \ref{thm:main}. The complete determination of ${\rm sp}_0(G)$ is reduced to a thorough understanding of $\lambda$. 

\begin{nota}
\rm
Let $\mathbb{N}=\{h\in\mathbb{Z}\colon h\ge 0\}$.
Let $p>2$ be a fixed prime number.

For a positive integer $k$, denote $\mathbb{Z}/p^k\mathbb{Z}$ by $\mathbb{Z}_{p^k}$; it is a quotient ring of $\mathbb{Z}$, and further is a field when $k=1$.

For an integer $a$, let $\delta_{a}=1$ (resp. $\delta_a=0$) if $a$ is even (resp. odd).

Given subsets $U_i\subseteq\mathbb{Q}$ and $a_i\in\mathbb{Q}$, $i=1,\ldots,r$, let $\sum_{i=1}^ra_iU_i\subseteq\mathbb{Q}$ denote the set
$\big\{\sum_{i=1}^ra_iu_i\colon u_i\in U_i, 1\le i\le r\big\}.$
\end{nota}

Recall the basic fact that each nontrivial finite abelian $p$-group is isomorphic to $\prod_{i=1}^{n}\mathbb{Z}_{p^{k_i}}$ for a unique tuple $(k_{1},\ldots,k_{n})$ with $k_1\geq\cdots\geq k_n\geq 1$.

\begin{defn}\label{defn:type}
\rm Call a tuple of integers $(k_1,\ldots,k_n)$ with $k_1\ge\cdots\ge k_n\ge 0$ a {\it generalized type}. When $k_{i+1}=\cdots=k_{i+s}=q$, sometimes we abbreviate the segment $k_{i+1},\ldots,k_{i+s}$ to $q^{[s]}$.

Given generalized types $\alpha=(k_{1},\ldots,k_{n})$ and $\beta=(\ell_{1},\ldots,\ell_{m})$, denote $\beta\leq\alpha$ if $m\leq n$ and $\ell_{i}\leq k_{i}$ for all $i\in\{1,\ldots,m\}$; denote $\beta<\alpha$ if $\beta\leq\alpha$ and $\beta\ne\alpha$.

A generalized type $(k_1,\ldots,k_n)$ is called a {\it type} if $k_n\ge 1$.
\end{defn}

\begin{defn}
\rm Given a type $\alpha=(k_{1},\ldots,k_{n})$, let $e(\alpha)=k_{1}$ and $A(\alpha)={\prod}_{i=1}^{n}\mathbb{Z}_{p^{k_{i}}}$, so that the exponent of $A(\alpha)$ is $p^{e(\alpha)}$.

Denote ${\rm sp}(A(\alpha))$ as ${\rm sp}(\alpha)$, denote $\mu_0(A(\alpha))$ as $\mu_0(\alpha)$, and so forth.
\end{defn}

\medskip

Suppose $\tau=\big(q_1^{[s_1]},\ldots,q_c^{[s_c]}\big)$, with $q_1>\cdots>q_c\ge 1$ and each $s_i>0$.
Let $n=s_1+\cdots+s_c$, and $e=e(\tau)=q_1$.

For each $q_c\le z\le e$, let
\begin{align*}
b_z&=\min\{b\colon q_b\le z, \ 1\le b\le c\},  \\
\theta_z&=\big((z-q_c)^{[s_c]},\ldots,(z-q_{b_z})^{[s_{b_z}]}\big),  \\
\Omega_z&=\Big\{1+p^z\delta_{n-r}+{\sum}_{i=1}^rp^{d'_i}\colon (d'_1,\ldots,d'_r)\le\theta_z,\ d'_r=0, \ 1\le r\le n-n_{b_z-1}\Big\}.
\end{align*}
Put
\begin{align*}
\Gamma_z&=\Big\{{\sum}_{i=0}^{z-1}a_i(p^z-p^i)\colon 0\le a_i<p\Big\},  \qquad 1\le z\le e; \\
\Lambda_z&=\begin{cases}
\{-p^e(1+\delta_n)\}, &z=0, \\
p^{e-z}\Gamma_z+\{p^e(1-\delta_n)-2p^{e-z}\}, &1\le z<q_c, \\
p^{e-z}(\Gamma_z-\Omega_z), &q_c\le z\le e.\end{cases}
\end{align*}
For each $x$ with $0\le x<p^e$, set
\begin{align*}
\lambda(x)&=\min\Big\{h\in{\bigcup}_{z=0}^{e}\Lambda_z\colon h\equiv -x\ (\bmod{\ p^e})\Big\},  \\
\overline{\lambda}(x)&=\frac{1}{2}(np^e+\lambda(x)).
\end{align*}

With the above notations, the main result of the paper is stated as follows.
\begin{thm}\label{thm:main}
{\rm(a)} The reduced genus spectrum is given by
\begin{align*}
{\rm sp}_0(\tau)=p^e\mathbb{N}+\big\{\frac{n}{2}p^e\big\}+\frac{1}{2}{\bigcup}_{z=0}^{e}\Lambda_z
=p^e\mathbb{N}+\big\{\overline{\lambda}(x)\colon 0\le x<p^e\big\}.
\end{align*}
{\rm(b)} The reduced minimum genus $\mu_0(\tau)$ equals
$$\frac{n-1}{2}p^e-\frac{1}{2}\max\Big\{\delta_np^e,\max\Big\{p^{e-q_b}+{\sum}_{v=b}^cs_vp^{e-q_v}\colon 1\le b\le c,\
2\mid n_{b-1}\Big\}\Big\}.$$
{\rm(c)} Let $\lambda_\ast=\max_{0\le x<p^e}\lambda(x)$, and $\overline{\lambda}_\ast=(np^e+\lambda_\ast)/2$.
Then the reduced stable upper genus is given by
$$\sigma_0(\tau)=\min\big\{w\in{\rm sp}_0(\tau)\colon w>\overline{\lambda}_\ast-p^e\big\}
=\min_{0\le x<p^e}\Big(p^e\Big[\frac{\overline{\lambda}_\ast-\overline{\lambda}(x)}{p^e}\Big]+\overline{\lambda}(x)\Big).$$
\end{thm}

\begin{rmk}
\rm
It is worth highlighting that each $\Lambda_z$ is a finite set, so that each $\lambda(x)$ can be found in finitely many steps. The complexity of ${\rm sp}_0(\tau)$ is captured by the finite set $\big\{\overline{\lambda}(x)\colon 0\le x<p^e\big\}$.

The formula in (b) generalizes the one for $\mu_0$ given in \cite{MS19} Theorem 6.2, and also, similarly as claimed in \cite{MS19} Section 5.7, improves \cite{Ma65} Theorem 4 in the case of abelian $p$-groups.
\end{rmk}

\begin{cor}\label{cor:no-gap}
A sufficient and necessary condition for $\mu_0(\tau)=\sigma_0(\tau)$ is that $\lambda(x)>\lambda_\ast-2p^e$ for all $0\le x<p^e$, or equivalently, the set $\{\lambda(x)\colon 0\le x<p^e\}$ is contained in an interval of length $2p^e$.
\end{cor}

Expecting a complete solution to the Hurwitz Problem for abelian $p$-groups, we leave the following problem to future research:
\begin{prob}
Find a useful formula for $\lambda(x)$.
\end{prob}

In Section \ref{sec:group} we develop some group-theoretic techniques, to give a precise expression for elements of the genus spectrum.
In Section \ref{sec:combinatoric}, we carefully look into the combinatorial structure of the genus spectrum, and prove Theorem \ref{thm:main}. In Section \ref{sec:example} we reproduce some existing results without much effort, so as to illustrate the power and verify the correctness of Theorem \ref{thm:main}.

\section{Group-theoretic aspect} \label{sec:group}

Suppose $K$ is a finite abelian $p$-group.
Let ${\rm rk}(K)$ denote the minimal number of generators of $K$. For $s\in\mathbb{N}$, let $p^sK=\{p^s\mathfrak{x}\colon \mathfrak{x}\in K\}$.
For $\mathfrak{x}\in K$, let $|\mathfrak{x}|$ denote its order.
If $K\cong A(\alpha)$, then we say that $K$ is of {\it type $\alpha$} and write ${\rm tp}(K)=\alpha$.

\begin{lem} \label{lem:basic}
If $A$ is a finite abelian $p$-group and $B$ is a proper subgroup of $A$, then ${\rm tp}(B)<{\rm tp}(A)$ and ${\rm tp}(A/B)<{\rm tp}(A)$.
\end{lem}

Although this should be well-known, we give a proof for completeness.
\begin{proof}

Suppose ${\rm tp}(A)=(k_1,\ldots,k_n)$ and ${\rm tp}(B)=(\ell_{1},\ldots,\ell_{m})$.

Observe that $A_\diamond:=\{\mathfrak{a}\in A\colon p\mathfrak{a}=0\}\cong\mathbb{Z}_p^n$, and $B_\diamond:=\{\mathfrak{b}\in B\colon p\mathfrak{b}=0\}\cong\mathbb{Z}_p^m$. Clearly, $B_\diamond\le A_\diamond$, so $m\le n$, i.e. ${\rm rk}(B)\le{\rm rk}(A)$.

For a finite set $Y$, let $\#Y$ denote its cardinality.

For each $s\ge 1$, ${\rm rk}(p^sA)=\#\{i\colon k_i>s\}$, ${\rm rk}(p^sB)=\#\{i\colon \ell_i>s\}$. Since $p^sB\le p^sA$, by the previous paragraph, we have ${\rm rk}(p^sB)\le{\rm rk}(p^sA)$. Hence $\#\{i\colon \ell_i>s\}\le\#\{i\colon k_i>s\}$ for all $s$, which combined with the assumption that $B$ is proper implies ${\rm tp}(B)<{\rm tp}(A)$.

Observe that
$$A^\ast:=\hom(A,\mathbb{Z}_{p^{k_1}})\cong A, \qquad (A/B)^\ast:=\hom(A/B,\mathbb{Z}_{p^{k_1}})\cong A/B,$$
and the quotient map $A\twoheadrightarrow A/B$ induces an embedding $(A/B)^\ast\hookrightarrow A^\ast$.
By the above, ${\rm tp}((A/B)^\ast)<{\rm tp}(A^\ast)$, which is the same as ${\rm tp}(A/B)<{\rm tp}(A)$.
\end{proof}

From now on, let $G$ denote a finite abelian $p$-group of type $\tau$.

Refer to \cite{MS19} Page 5 or \cite{Ta11} Page 3 for the following fundamental result:
\begin{lem}\label{lem:fundamental}
An integer $g$ belongs to ${\rm sp}(\tau)$ if and only if there exist elements $\mathfrak{x}_1,\mathfrak{y}_1,\ldots,\mathfrak{x}_b,\mathfrak{y}_b,\mathfrak{z}_1,\ldots,\mathfrak{z}_{s}\in G$ such that
\begin{align}
\langle \mathfrak{x}_1,\mathfrak{y}_1,\ldots,\mathfrak{x}_b,\mathfrak{y}_b,\mathfrak{z}_1,\ldots,\mathfrak{z}_{s}\rangle&=G, \label{eq:generate} \\
\mathfrak{z}_1+\cdots+\mathfrak{z}_{s}&=0, \label{eq:sum-to-zero} \\
b-1+\frac{1}{2}\sum\limits_{i=1}^{s}\Big(1-\frac{1}{|\mathfrak{z}_{i}|}\Big)&=\frac{g-1}{|G|}.  \label{eq:well-known}
\end{align}
\end{lem}

\medskip

In the notation of Lemma \ref{lem:fundamental}, let $N=\langle \mathfrak{z}_1,\ldots,\mathfrak{z}_{s}\rangle$.

If $s=0$, then $N=0$, and (\ref{eq:generate}) holds for some $\mathfrak{x}_1,\mathfrak{y}_1,\ldots,\mathfrak{x}_b,\mathfrak{y}_b$ if and only if $2b\ge{\rm rk}(G)$. So
\begin{align}
\frac{g-1}{|G|}=b-1\ge\Big[\frac{{\rm rk}(G)-1}{2}\Big].  \label{eq:s=0}
\end{align}

Since $s\ne 1$, now we assume $s\ge 2$ and $N\ne 0$. Suppose ${\rm tp}(N)=\alpha$.
In virtue of (\ref{eq:sum-to-zero}), there are at least two $i$'s with $|\mathfrak{z}_i|=p^{e(\alpha)}$; without loss of generality we may assume $|\mathfrak{z}_{s}|=|\mathfrak{z}_{s-1}|=p^{e(\alpha)}$. Clearly, $N=\langle\mathfrak{z}_1,\ldots,\mathfrak{z}_{s-1}\rangle$.
Let $\{\mathfrak{z}_{v_1},\ldots,\mathfrak{z}_{v_r}\}\subseteq\{\mathfrak{z}_1,\ldots,\mathfrak{z}_{s-1}\}$ be a minimal generating set for $N$,
with $|\mathfrak{z}_{v_i}|=p^{d_i}$.
Then $\sum_{i=1}^{s}(1-1/|\mathfrak{z}_i|)$ in (\ref{eq:well-known}) can be rewritten as the sum of
$$1-\frac{1}{p^{e(\alpha)}}, \qquad  \sum_{i=1}^{r}\Big(1-\frac{1}{p^{d_i}}\Big), \qquad  \sum_{t=1}^{s-1-r}\Big(1-\frac{1}{p^{c_t}}\Big),$$
which are respectively contributed by $\mathfrak{z}_{s}$, the elements in the generating set $\{\mathfrak{z}_{v_1},\ldots,\mathfrak{z}_{v_r}\}$, and those in $\{\mathfrak{z}_1,\ldots,\mathfrak{z}_{s-1}\}\setminus\{\mathfrak{z}_{v_1},\ldots,\mathfrak{z}_{v_r}\}$.

Obviously, (\ref{eq:generate}) is equivalent to $\langle \overline{\mathfrak{x}_1},\overline{\mathfrak{y}_1},\ldots,\overline{\mathfrak{x}_b},\overline{\mathfrak{y}_b}\rangle=G/N$, where $\overline{\mathfrak{a}}$ denotes the image of $\mathfrak{a}\in G$ under the quotient map $G\twoheadrightarrow G/N$. Let $u$ denote the minimum value of ${\rm rk}(G/N)$ when $N$ ranges over all subgroups of $G$ with ${\rm tp}(N)=\alpha$. Then there exist $\mathfrak{x}_1,\mathfrak{y}_1,\ldots,\mathfrak{x}_b,\mathfrak{y}_b$ fulfilling (\ref{eq:generate}) if and only if $2b\ge u$, i.e. $b=[(u+1)/2]+h$ for some $h\in\mathbb{N}$.

Therefore, the Hurwitz Problem is reduced to
\begin{prob} \label{prob:key}
\rm For each type $\alpha\leq \tau$,
\begin{enumerate}
  \item[\rm (i)] determine all possible tuples $(d_{1},\ldots,d_{r})$ such that there is a minimal generating set
                 $\{\mathfrak{z}_{1},\ldots,\mathfrak{z}_{r}\}$ for $A(\alpha)$ with $|\mathfrak{z}_{i}|=p^{d_{i}}$;
  \item[\rm (ii)] determine $u=\min\{{\rm rk}(A(\tau)/N)\colon N\le A(\tau), \ {\rm tp}(N)=\alpha\}$.
\end{enumerate}
\end{prob}

Once this problem is solved, the integers $g\in{\rm sp}(\tau)$ other than those given by (\ref{eq:s=0}) can be found via
\begin{align}
\frac{2(g-1)}{|G|}=2h+u-\delta_u-\frac{1}{p^{e(\alpha)}}+\sum\limits_{i=1}^{r}\Big(1-\frac{1}{p^{d_{i}}}\Big)
+\sum\limits_{t\in\Omega}\Big(1-\frac{1}{p^{c_t}}\Big)  \label{eq:equation'}
\end{align}
for varying $\alpha$, $(d_{1},\ldots,d_{r})$, $h\in\mathbb{N}$, finite sets $\Omega$ (which may be empty), and $c_t\in\{1,\ldots,e(\alpha)\}$.
We can alternatively write (\ref{eq:equation'}) as
\begin{align}
\frac{2(g-1)}{|G|}=2h'+u-\delta_u-\frac{1}{p^{e(\alpha)}}+\sum\limits_{i=1}^{r}\Big(1-\frac{1}{p^{d_{i}}}\Big)
+\sum\limits_{j=1}^{e(\alpha)}a_j\Big(1-\frac{1}{p^{j}}\Big)  \label{eq:equation}
\end{align}
with $h'\in\mathbb{N}$ and $a_j\in\{0,1,\ldots,p-1\}$.


The two parts of Problem \ref{prob:key} are separately settled in the following two theorems.

\begin{thm}\label{thm1}
Let $\alpha=(\ell_{1},\ldots,\ell_{m})$ be a type.
\begin{enumerate}
  \item [\rm (a)] If $\{\mathfrak{z}_1,\ldots,\mathfrak{z}_r\}$ is a minimal generating set for $A(\alpha)$ with $|\mathfrak{z}_{i}|=p^{d_{i}}$ and
        $\ell_1=d_1\ge\cdots\ge d_r\ge 1$, then $r=m$ and $d_{i}\ge\ell_{i}$ for all $i$.
  \item [\rm (b)] Conversely, given $d_{1},\ldots,d_{m}$ such that $\ell_1=d_1\ge\cdots\ge d_m\ge 1$ and $d_{i}\ge\ell_{i}$ for all $i$,
        there exists a minimal generating set $\{\mathfrak{z}_{1},\ldots,\mathfrak{z}_{m}\}$ for $A(\alpha)$ with $|\mathfrak{z}_{i}|=p^{d_{i}}$.
\end{enumerate}
\end{thm}

\begin{proof}
(a) Consider the reduction map $\mathcal{R}:A(\alpha)\twoheadrightarrow A(\alpha)/pA(\alpha)\cong\mathbb{Z}_p^{m}.$
Since $\langle\mathcal{R}(\mathfrak{z}_1),\ldots,\mathcal{R}(\mathfrak{z}_r)\rangle=\mathbb{Z}_p^{m}$, we have $r\geq m$.

On the other hand, assume $r>m$. Then for some $j\in\{1,\ldots,r\}$, $\mathcal{R}(\mathfrak{z}_j)$ can be written as a linear combination of $\mathcal{R}(\mathfrak{z}_i), i\in\{1,\ldots,r\}-\{j\}$; without loss of generality we may assume $j=r$.
Suppose
$\mathcal{R}(\mathfrak{z}_{r})=\sum_{i=1}^{r-1}a_{i}\mathcal{R}(\mathfrak{z}_{i})$.
It follows that $\mathfrak{z}_{r}=\sum_{i=1}^{r-1}a_{i}\mathfrak{z}_{i}+\mathfrak{x}$ for some $\mathfrak{x}\in\ker(\mathcal{R})$. Hence $\mathfrak{x}=p\mathfrak{y}$ for some $\mathfrak{y}\in A(\alpha)$.
Writing $\mathfrak{y}=\sum_{i=1}^{r}c_{i}\mathfrak{z}_{i}$, we obtain
$$(1-pc_{r})\mathfrak{z}_{r}={\sum}_{i=1}^{r-1}(a_{i}+pc_{i})\mathfrak{z}_{i},$$
which, due to $\mathfrak{z}_r=\sum_{j=0}^{\ell_1-1}(pc_r)^j(1-pc_r)\mathfrak{z}_r$,
implies $\mathfrak{z}_{r}\in\langle \mathfrak{z}_{1},\ldots,\mathfrak{z}_{r-1}\rangle$.
This contradicts the minimality of $\{\mathfrak{z}_{1},\ldots,\mathfrak{z}_{r}\}$. Thus $r=m$.

It follows from $A(\alpha)=\langle \mathfrak{z}_1,\ldots,\mathfrak{z}_m\rangle$ that $A(\alpha)$ is a quotient of $\prod_{i=1}^m\mathbb{Z}_{p^{d_i}}$. Hence by Lemma \ref{lem:basic}, $\ell_i\le d_i$ for all $i$.

(b) Suppose $\ell_1=d_1\ge\cdots\ge d_m\ge 1$ and $d_{i}\ge\ell_{i}$ for all $i$.

Let $i_0=\max\{i\colon\ell_i=\ell_1\}$. Take
$$\mathfrak{z}_i=\begin{cases} \mathfrak{f}_i, &1\le i\le i_0, \\ p^{\ell_1-d_i}\mathfrak{f}_1+\mathfrak{f}_i, &i_0<i\le m, \end{cases}$$
where $\mathfrak{f}_i\in A(\alpha)$ is the element with the $i$-th entry being $1$ and the others being $0$.
Then $\{\mathfrak{z}_1,\ldots,\mathfrak{z}_m\}$ is a minimal generating set, with $|\mathfrak{z}_i|=p^{d_i}$.
\end{proof}

\begin{thm}\label{thm2}
Suppose $\alpha=(\ell_{1},\ldots,\ell_{m})\le\tau=(k_{1},\ldots,k_{n})$. Let
$$u^\tau_\alpha=\max\big\{j-f^\tau_\alpha(j)\colon 1\leq j\leq n\big\},$$
where $f^\tau_\alpha(j)=\max\{i\colon \ell_{i}\geq k_{j}\ \text{or\ }i=0\}$.
Then
$$\min\big\{{\rm rk}(A(\tau)/N)\colon N\le A(\tau),\ {\rm tp}(N)=\alpha\big\}=u^\tau_\alpha.$$
\end{thm}

\begin{proof}
Let $k=k_1$. Embed $A(\tau)$ into $\mathbb{Z}_{p^k}^{n}$ via
\begin{align*}
\Phi:A(\tau)\hookrightarrow \mathbb{Z}_{p^k}^{n}, \ \ \ \
(c_1,\ldots,c_n)\mapsto (p^{k-k_{1}}c_1,\ldots,p^{k-k_{n}}c_n).
\end{align*}
The image ${\rm Im}(\Phi)$ consists of elements $(a_{1},\ldots,a_{n})\in\mathbb{Z}_{p^{k}}^{n}$ with
\begin{align}
p^{k-k_{j}}\mid a_{j}, \qquad 1\leq j\leq n.    \label{eq:elements}
\end{align}
We identify $A(\tau)$ with ${\rm Im}(\Phi)$ and write elements of $A(\tau)$ as row-vectors.

For a matrix $X$ over some ring $R$, let $X_i$ denote its $i$-th row, and let $X_{ij}$ or $X_{i,j}$ denote the $(i,j)$-entry; let ${\rm rk}(X)$ denote its rank if $R$ is a field.

(i) Suppose $N\le A(\tau)$ has type $\alpha$. We show ${\rm rk}(A(\tau)/N)\ge u^\tau_\alpha$.

Now that $N$ is a subgroup of $\mathbb{Z}_{p^{k}}^{n}$ of type $\alpha$, by \cite{CS13} Theorem 3.9, there exists $Q\in{\rm GL}(n,\mathbb{Z}_{p^{k}})$ such that
$$N=\langle p^{k-\ell_1}Q_1,\ldots,p^{k-\ell_m}Q_m\rangle.$$
Since $p^{k-\ell_i}Q_i\in A(\tau)$, by (\ref{eq:elements}) we have $p^{\ell_i-k_j}\mid Q_{i,j}$ whenever $\ell_i>k_j$.

Let $v={\rm rk}(A(\tau)/N)$. Then
\begin{align*}
\mathbb{Z}_p^{v}\cong\frac{A(\tau)/N}{p(A(\tau)/N)}=\frac{A(\tau)/N}{(pA(\tau)+N)/N}\cong\frac{A(\tau)}{pA(\tau)+N}.
\end{align*}
Since $A(\tau)/pA(\tau)\cong\mathbb{Z}_p^n$, we have
\begin{align}
\frac{N}{N\cap pA(\tau)}\cong\frac{pA(\tau)+N}{pA(\tau)}\cong\mathbb{Z}_p^{n-v}.   \label{eq:quotient}
\end{align}

Let $\Delta^\tau_\alpha=\{(i,j)\colon \ell_{i}\geq k_{j}\}$. Define $Q^\vee\in{\rm GL}(n,\mathbb{Z}_p)$ by setting
$$Q^\vee_{ij}=\begin{cases}  Q_{ij}/p^{\ell_{i}-k_{j}}, &(i,j)\in\Delta^\tau_\alpha, \\  0,&\text{otherwise}. \end{cases}$$
The map
$$\mathbb{Z}_p^n\to \frac{N}{N\cap pA(\tau)}, \qquad    \vec{c}=(c_1,\ldots,c_n)\mapsto \sum_{i=1}^mc_ip^{k-\ell_i}Q_i$$
is well-defined and surjective, whose kernel is $\{\vec{c}\in\mathbb{Z}_p^n\colon\vec{c}Q^\vee=0\}$.
Hence by (\ref{eq:quotient}), ${\rm rk}(Q^\vee)=n-v$.

For each $j$, since $\Delta^\tau_\alpha$ is contained in the disjoint union of the two sets
$$\{(i',j')\colon 1\le i'\le f^\tau_\alpha(j),\ 1\le j'\le j\}, \qquad \{(i',j')\colon j<j'\le n\},$$
which respectively have $f^\tau_\alpha(j)$ rows and $n-j$ columns, we have ${\rm rk}(Q^\vee)\le f^\tau_\alpha(j)+n-j$. Hence
$${\rm rk}(Q^\vee)\le\min\big\{n+f^\tau_\alpha(j)-j\colon 1\leq j\leq n\big\}=n-u^\tau_\alpha.$$
Thus ${\rm rk}(A(\tau)/N)=v\ge u^\tau_\alpha$.

(ii) Take $Q\in{\rm GL}(n,\mathbb{Z}_{p^k})$ with
$$Q_{ij}=\begin{cases}  1, &1\le i=j\le n, \\ p^{\ell_i-k_{j}}, &j=u^\tau_\alpha+i,\ 1\le i\le n-u^\tau_\alpha, \\ 0, &\text{otherwise}.\end{cases}$$
Then ${\rm rk}(Q^\vee)=n-u^\tau_\alpha$, so that $N=\langle p^{\ell_1}Q_1,\ldots,p^{\ell_m}Q_m\rangle$ realizes the lower bound:
${\rm rk}(A(\tau)/N)=u^\tau_\alpha$.
\end{proof}

\section{Combinatorial aspect} \label{sec:combinatoric}

Fix $\tau=(k_1,\ldots,k_n)=\big(q_1^{[s_1]},\ldots,q_c^{[s_c]}\big)$, with $q_1>\cdots>q_c\ge 1$ and $s_1+\cdots+s_c=n$. For each $v\in\{1,\ldots,c\}$, let $n_v=s_1+\cdots+s_v$; set $n_0=0$. Let $e=q_1=k_1$.

Recall the notations given before stating Theorem \ref{thm:main}.

Given $u\in\{0,\ldots,n\}$, $z\in\{1,\ldots,e\}$, let
$$\Xi_{u,z}=\{\text{types\ }\alpha\colon \alpha\le\tau,\ u^\tau_\alpha=u,\ e(\alpha)=z\}.$$
For $\alpha=(\ell_1,\ldots,\ell_m)\in\Xi_{u,z}$, let
$$U(\alpha)=2\mathbb{N}+\{u-\delta_u-p^{-z}\}+V(\alpha)+p^{-z}\Gamma_z,$$
where
\begin{align*}
V(\alpha)=\Big\{{\sum}_{i=1}^m(1-p^{-d_i})\colon \alpha\le (d_1,\ldots,d_m), \ d_1=z\Big\}.
\end{align*}
In view of (\ref{eq:s=0}), (\ref{eq:equation}), Theorem \ref{thm1}, Theorem \ref{thm2},
the determination of the genus spectrum is reduced to
\begin{align}
{\rm sp}_0(\tau)=\Big(p^e\mathbb{N}+\big\{\frac{n}{2}p^e\big\}+\frac{1}{2}\Lambda_0\Big)\cup\frac{1}{2}p^eS,   \qquad
S={\bigcup}_{u=0}^n{\bigcup}_{z=1}^{e}{\bigcup}_{\alpha\in\Xi_{u,z}}U(\alpha).   \label{eq:sp}
\end{align}

Observe that if $\alpha,\alpha'\in\Xi_{u,z}$ and $\alpha'\le\alpha$, then $V(\alpha)\subseteq V(\alpha')$. It follows that if $\nu$ is the minimal element of $\Xi_{u,z}$, then ${\bigcup}_{\alpha\in\Xi_{u,z}}U(\alpha)=U(\nu)$. Fortunately, the minimal element does exist and is easy to find.

\begin{lem}\label{lem:min}
{\rm(a)} $\Xi_{n,z}\ne\emptyset$ if and only if $z<q_c$, in which case $(z)$ is the minimal element of $\Xi_{n,z}$.

{\rm(b)} Suppose $n_{b-1}\le u<n_b$, with $1\le b\le c$. Then $\Xi_{u,z}\ne\emptyset$ if and only if $q_b\le z\le e$, in which case
$$\nu_{u,z}:=\big(z,q_b^{[n_b-u-1]},q_{b+1}^{[s_{b+1}]},\ldots,q_c^{[s_c]}\big)$$
is the minimal element of $\Xi_{u,z}$.
\end{lem}

\begin{proof}
(a) Recalling the definitions of $u^\tau_\alpha$ and $f^\tau_\alpha$ (in Theorem \ref{thm2}), we see that $u^\tau_\alpha=n$ is equivalent to $f^\tau_\alpha(n)=0$, which in turn is equivalent to $e(\alpha)<q_c$.

If $\alpha\in\Xi_{n,z}$, i.e. $u^\tau_\alpha=n$ and $e(\alpha)=z$, then $z<q_c$. On the other hand, whenever $z<q_c$, we have $(z)\in\Xi_{n,z}$, and $(z)\le\alpha$ for each $\alpha\in\Xi_{n,z}$.

(b) Suppose $\alpha=(\ell_1,\ldots,\ell_m)\in\Xi_{u,z}$. Since $u^\tau_\alpha=u$, we have $f^\tau_\alpha(j)\ge j-u$ for all $j$.
In particular, $f^\tau_\alpha(u+1)\ge 1$, implying $z=\ell_1\ge k_{u+1}=q_b$.

Now suppose $q_b\le z\le e$. If $n_{t-1}<j\le n_t$ with $b\le t\le c$, then $k_j=q_t$, so that $f^\tau_{\nu_{u,z}}(j)=n_b-u+(s_{b+1}+\cdots+s_t)=n_t-u$.
If $j\le n_{b-1}$, then $j-f^\tau_{\nu_{u,z}}(j)\le n_{b-1}\le u$.
Hence
$$u^\tau_{\nu_{u,z}}=\max\{j-f^\tau_{\nu_{u,z}}(j)\colon 1\le j\le n\}=u,$$
showing $\nu_{u,z}\in\Xi_{u,z}$.

Given $\alpha=(\ell_1,\ldots,\ell_m)\in\Xi_{u,z}$, to show $\nu_{u,z}\le\alpha$, we put $n'_{b-1}=1$ and $n'_t=n_t-u$ for $b\le t\le c$. There are $n'_c$ entries in $\nu_{u,z}$.
For each $i$ with $n'_{t-1}<i\le n'_t$, $t\in\{b,\ldots,c\}$, the $i$-th entry of $\nu_{u,z}$ is $q_t$; since $f^\tau_\alpha(n_t)\ge n'_t$, one has $\ell_{n'_t}\ge k_{n_t}=q_t$, so that $\ell_i\ge \ell_{n'_t}\ge q_t$. Thus, $\nu_{u,z}\le\alpha$.
\end{proof}

\medskip

\begin{proof}[Proof of Theorem \ref{thm:main}]
(a) From Lemma \ref{lem:min}, it follows that $S=S_1\cup S_2$, where
\begin{align}
S_1&={\bigcup}_{z=1}^{q_c-1}U((z)),  \label{eq:S1} \\
S_2&={\bigcup}_{u=0}^{n-1}{\bigcup}_{z=1}^{e}U(\nu_{u,z})
={\bigcup}_{b=1}^c{\bigcup}_{u=n_{b-1}}^{n_b-1}{\bigcup}_{z=q_b}^{e}U(\nu_{u,z})  \nonumber  \\
&={\bigcup}_{z=q_c}^{e}{\bigcup}_{b=b_z}^{c}{\bigcup}_{u=n_{b-1}}^{n_b-1}U(\nu_{u,z})
={\bigcup}_{z=q_c}^{e}{\bigcup}_{u=n_{b_z-1}}^{n-1}U(\nu_{u,z}).  \label{eq:S2}
\end{align}

Suppose $q_c\le z\le e$.
We can rewrite
${\sum}_{i=1}^{n-u}(1-p^{-d_i})$ as $n-u-p^{-z}{\sum}_{i=1}^{n-u}p^{d'_i}$, with $d'_i=z-d_{n-u+1-i}$.
Let
$$\overline{\nu}_{u,z}=\big((z-q_c)^{[s_c]},\ldots,(z-q_{b_z+1})^{[s_{b_z+1}]},(z-q_{b_z})^{[n_{b_z}-u]}\big).$$
Then $\nu_{u,z}\le(d_1,\ldots,d_{n-u})$, $d_1=z$ if and only if $(d'_1,\ldots,d'_{n-u})\le\overline{\nu}_{u,z}$, $d'_{n-u}=0$.
(Here $(d'_1,\ldots,d'_{n-u})$ and $\overline{\nu}_{u,z}$ are generalized types,
so $d'_i$'s as well as entries of $\overline{\nu}_{u,z}$ are allowed to vanish.)
Thus,
$$V(\nu_{u,z})=\{n-u\}-p^{-z}\Big\{{\sum}_{i=1}^{n-u}p^{d'_i}\colon (d'_1,\ldots,d'_{n-u})\le\overline{\nu}_{u,z},\ d'_{n-u}=0\Big\}.$$
Consequently,
\begin{align*}
{\bigcup}_{u=n_{b_z-1}}^{n-1}U(\nu_{u,z})&=2\mathbb{N}+{\bigcup}_{u=n_{b_z-1}}^{n-1}\big(\{u-\delta_u-p^{-z}\}+V(\nu_{u,z})\big)+p^{-z}\Gamma_z   \\
&=2\mathbb{N}+\{n\}+p^{-z}(\Gamma_z-\Omega_z).
\end{align*}

As a complement, for $1\le z<q_c$, one can check that
$$U((z))=2\mathbb{N}+\{n+1-\delta_n-2p^{-z}\}+p^{-z}\Gamma_z.$$

Combing (\ref{eq:sp}), (\ref{eq:S1}), (\ref{eq:S2}), we obtain
$${\rm sp}_0(\tau)=p^e\mathbb{N}+\big\{\frac{n}{2}p^e\big\}+\frac{1}{2}{\bigcup}_{z=0}^{e}\Lambda_z.$$

To show the second equality, let
$$M_x=\Big\{h\in{\bigcup}_{z=0}^{e}\Lambda_z\colon h\equiv -x\ (\bmod{\ p^e})\Big\}.$$
Observe that $M_x\subset2p^e\mathbb{N}+\{\lambda(x)\}$.
Since
$${\bigcup}_{z=0}^{e}\Lambda_z={\bigcup}_{x=0}^{p^e-1}M_x\subset 2p^e\mathbb{N}+\{\lambda(x)\colon 0\le x<p^e\},$$
we have
$$p^e\mathbb{N}+\big\{\frac{n}{2}p^e\big\}+\frac{1}{2}{\bigcup}_{z=0}^{e}\Lambda_z
\subseteq p^e\mathbb{N}+\big\{\overline{\lambda}(x)\colon 0\le x<p^e\big\}.$$
The other direction $\supseteq$ is straightforward.

\medskip

(b) For each $1\le z<q_c$,
\begin{align*}
\min\Lambda_z&=p^{e-z}\min\Gamma_z+p^e(1-\delta_n)-2p^{e-z}  \\
&\ge p^e(1-\delta_n)-2p^{e-z}>-p^e(1+\delta_n).
\end{align*}
Hence $\min_{0\le z<q_c}\min\Lambda_z=-p^e(1+\delta_n).$

For each $b$ with $2\mid n_{b-1}$, set
$$w(b)=p^{e-q_b}+{\sum}_{v=b}^cs_vp^{e-q_v}.$$

Suppose $q_c\le z\le q_1$. Let $y\in p^{e-z}\Omega_z$ be an arbitrary element, say
$$y=p^{e-z}+p^e\delta_{n-r}+p^{e-z}{\sum}_{i=1}^rp^{d'_i},$$
with $(d'_1,\ldots,d'_r)\le\theta_z$ and $d'_r=0$. 

If $2\mid n_{b_z-1}$, then $y\le p^e+w(b_z).$

If $2\nmid n_{b_z-1}$, then either $n\not\equiv r\pmod{2}$ or $(d'_1,\ldots,d'_r)<\theta_z$; each implies
        $y<1+p^e+{\sum}_{v=1}^cs_vq^{e-q_v}=p^e+w(1).$

Hence
$${\max}_{q_c\le z\le e}\max p^{e-z}\Omega_z=p^e+\max\big\{w(b)\colon 1\le b\le c,\ 2\mid n_{b-1}\big\}.$$
As $\min \Gamma_z=0$ for each $z$, the assertion follows from (a).

\medskip

(c) Let $w_0$ denote the minimal element of ${\rm sp}_0(\tau)$ larger than $\overline{\lambda}_\ast-p^e$.

By definition, $\overline{\lambda}_\ast-p^e\notin{\rm sp}_0(\tau)$, so $\overline{\lambda}_\ast-p^e<\sigma_0(\tau)$, i.e. $\sigma_0(\tau)\ge w_0$.

On the other hand, given $h\ge w_0$, there exists a unique $x\in\{0,\ldots,p^e-1\}$ such that $h\equiv\overline{\lambda}(x)\pmod{p^e}$.
Then
$h>\overline{\lambda}_\ast-p^e\ge\overline{\lambda}(x)-p^e$, implying
$h\in p^e\mathbb{N}+\overline{\lambda}(x)\subset{\rm sp}_0(\tau).$
This shows $\sigma_0(\tau)\le w_0$.

Therefore, $\sigma_0(\tau)=w_0$.
\end{proof}

\section{Recovering some existing results} \label{sec:example}

\begin{exmp}[Elementary $p$-groups]
\rm When $\tau=(1^{[n]})$, it is easy to see
$$\Lambda_0=\{-(\delta_n+1)p\}, \qquad \Lambda_1=\big\{a(p-1)-1-r-\delta_{n-r}p\colon 0\le a<p,\ 1\le r\le n\big\}.$$

\begin{itemize}
  \item If $n\ge p-1$, then writing $n+1=kp+d$, with $0\le d<p$, we have
        $$\lambda(x)=\begin{cases} -(k+\delta_{x+d})p-x,&0\le x\le d, \\ -(k-1+\delta_{x+d+1})p-x, &d<x<p.\end{cases}$$
        Clearly, $\lambda(d)\le\lambda(x)<\lambda(d)+2p$ for all $x$. Hence
        $$\mu_0(\tau)=\sigma_0(\tau)=\frac{1}{2}(n+1)(p-1)-p.$$
  \item If $n<p-1$, then $\lambda(0)=-(1+\delta_n)p$, $\lambda(1)=(p-2-n)p-1$, and
        $$\lambda(x)=\begin{cases} -\delta_{n-x+1}p-x, & 2\le x\le n+1, \\ (x-2-n)p-x,&n+2\le x<p.\end{cases}$$
        So $\lambda_\ast=\lambda(1)=(p-2-n)p-1$.
        Since
        $$\lambda(p-1)=(p-4-n)p+1>\lambda_\ast-2p,$$
        and $\lambda(p-1)\le\lambda(x)+2hp$ whenever $\lambda(x)+2hp>\lambda_\ast-2p$, by Theorem \ref{thm:main} (c) we have $\sigma_0(\tau)=\overline{\lambda}(p-1)$, i.e.
        $$\sigma_0(\tau)=\frac{1}{2}(p^2-4p+1).$$
        By Theorem \ref{thm:main} (b),
        $$\mu_0(\tau)=\frac{1}{2}\big((n-1)p-\max\{\delta_np,n+1\}\big),$$
        which is consistent with that $\min\{\lambda(0),\lambda(n+1)\}\le\lambda(x)$ for all $x$.
\end{itemize}

These recover \cite{MT98} Corollary 7.3. As stated in the end of \cite{MS19} Section 9.1, \cite{MT98} Corollary 7.3 (2) was at stake, due to the erroneous Remark in Section 7. Now we have confirmed it.
\end{exmp}

\begin{exmp}[Cyclic groups]
\rm Let $\tau=(e)$. We have $n=1$, $\theta_e=(0)$,
\begin{align*}
\Lambda_0&=\{-p^e\},  \\
\Lambda_z&=\Big\{{\sum}_{i=e-z}^{e-1}a_i(p^e-p^{i})+p^e-2p^{e-z}\colon 0\le a_i<p\Big\}, \qquad 0<z<e;  \\
\Lambda_e&=\Big\{{\sum}_{i=1}^{e-1}a_i(p^e-p^i)-p^e-2\colon 0\le a_i<p\Big\}.
\end{align*}

Since each $h\in\bigcup_{z=1}^e\Lambda_z$ with $h\equiv 0\pmod{p^e}$ is larger than $-p^e$, we have $\lambda(0)=-p^e$. For $0\le z<e$, each element of $\Lambda_z$ is a multiple of $p$; on the other hand, there is exactly one $h\in\Lambda_e$ satisfying $h\equiv -1\pmod{p^e}$, which is given by $a_0=\cdots=a_{e-1}=p-1$. Hence $\lambda(1)=(e(p-1)-2)p^e-1$.
Note that already $\lambda(1)=\max\bigcup_{z=0}^e\Lambda_z$, so $\lambda(1)\ge\lambda(x)$ for all $x$. Thus $\lambda_\ast=\lambda(1)$.

For $2\le x<p^e$, write $x=p^t\sum_{i=0}^{e-t-1}c_ip^i$ such that $c_0>0$ and $0\le c_i<p$ for all $i$. Let $f(x)=\sum_{i=0}^{e-t-1}c_i$.
Write $\lambda(x)=\lambda'(x)p^e-x$. 
\begin{itemize}
  \item If $t=0$ and $c_0\ge 2$, then $\lambda'(x)=f(x)-3$;
  \item if $t>0$ and $c_0\ge 2$, then $\lambda'(x)=f(x)-1$, realized by taking $a_{t+i}=c_i$ for $1\le i<e-t$ and $a_t=c_0-2$ in $\Lambda_{e-t}$;
  \item if $t=1$ and $c_0=1$, then $\lambda'(x)=f(x)+p-4$, realized by taking $a_i=c_{i-1}$ for $2\le i<e$, $a_1=0$ and $a_0=p-2$ in $\Lambda_e$;
  \item if $t>1$ and $c_0=1$, then $\lambda'(x)=f(x)+p-2$, realized by taking $a_{t+i}=c_i$ for $1\le i<e-t$, $a_t=0$ and $a_{t-1}=p-2$ in $\Lambda_{e-t+1}$.
\end{itemize}
In particular,
$$\lambda(p^e-1)=(e(p-1)-4)p^e+1>(e(p-1)-4)p^e-1=\lambda_\ast-2p^e.$$
The unique integer between $\lambda_\ast-2p^e$ and $\lambda(p^e-1)$ is $(e(p-1)-4)p^e$, but it does not belong to $\bigcup_{z=0}^e\Lambda_z$.
By Theorem \ref{thm:main} (c),
$$\sigma_0(\tau)=\overline{\lambda}(p^e-1)=\frac{1}{2}\big((e(p-1)-3)p^e+1\big),$$
which recovers \cite{KM91} Corollary 5.3.
\end{exmp}

\begin{exmp}[Groups with large invariants]
\rm Suppose $c=e$, $q_v=e+1-v$, ($1\le v\le e$), $s_1\ge p-2$, and $s_v\ge p-1$ for $2\le v\le e$. Then $n=s_1+\cdots+s_e$.
For each $1\le z\le e$, we have $b_z=e-z+1$, $\theta_z=\big((z-1)^{s_e},\ldots,0^{s_{e-z+1}}\big)$, and
\begin{align*}
\Lambda_z&=p^{e-z}\Big\{{\sum}_{i=0}^{z-1}a_i(p^z-p^i)\colon 0\le a_i<p\Big\}  \\
&\ \ \ \ -p^{e-z}\Big\{1+p^z\delta_{n-r}+{\sum}_{i=1}^rp^{d'_i}\colon (d'_1,\ldots,d'_r)\le\theta_z, d'_r=0,
1\le r\le {\sum}_{v=0}^{z-1}s_{e-v}\Big\}.
\end{align*}

Let $h=\sum_{v=1}^es_vp^{v-1}$. We are going to show that $|\lambda(x)+h+1|\le p^e$ for all $0\le x<p^e$.

For $1\le z\le e$, given $(d'_1,\ldots,d'_r)\le\theta_z$ with $d'_1=1$ and $1\le r\le n-n_{b_z-1}=n-(s_1+\cdots+s_{e-z})$, we can see that $$\rho:=p^{e-z}+p^e\delta_{n-r}+p^{e-z}{\sum}_{i=1}^rp^{d'_i}>1+p^e+h$$
only if $r=n-(s_1+\cdots+s_{e-z})$, $2\mid n-r$, and $\sum_{j=1}^{e-z}s_jp^{j-1}\le p^{e-z}-2$; but the last condition forces $s_2=\cdots=s_{e-z}=p-1$ and $s_1=p-2$, which contradicts $2\mid n-r$.
Hence $\rho\le 1+p^e+h.$
Consequently,
$$\min\Lambda_z\ge -p^{e-z}\max\Omega_z\ge-h-1-p^e.$$
Furthermore, note that $-p^e(1+\delta_n)\ge -h-1-p^e$, as $p^e\delta_n\le h+1$.
Therefore, $\lambda(x)\ge -h-1-p^e$ for all $x$.

It is sufficient to show $\Upsilon_x\ne\emptyset$ for each $x$, where
$$\Upsilon_x=\big\{w\in\Lambda_e\colon w\equiv -x\pmod{p^e}, \ w\le -h-1+p^e\big\}.$$

Given $0\le x<p^e$, let $x'$ be the remainder when dividing $h+1-x$ by $p^e$.
Write $x'=\sum_{i=1}^{e-1}c_ip^i$, with $0\le c_i<p$. Then $h-x'=\sum_{i=0}^{e-1}(s_{i+1}-c_i)p^i$.

If $s_1>c_0$, then $h-x'$ can be written as $\sum_{i=0}^rp^{d'_i}$ for some $(d'_1,\ldots,d'_r)\le\theta_e$ with $d'_r=0$, hence $x'-h-1-p^e\delta_{n-r}\in\Upsilon_x$, justifying $\Upsilon_x\ne\emptyset$.

Suppose $s_1\le c_0$, and suppose there exists $t$ such that $c_t<s_{t+1}$ and $c_i=s_{i+1}$ (equal to $p-1$) for $1\le i<t$.
\begin{itemize}
  \item If one of the two cases occurs: $c_0=p-1$; $s_1=c_0=p-2$ and $c_j>0$ for some $j\in\{t,\ldots,e-1\}$, then
        \begin{align*}
        h-x'
        &=(s_1-c_0+p)+{\sum}_{i=1}^{t-1}(p-1)p^i+(s_{t+1}-c_t-1)p^t \\
        &\ \ \ \ +{\sum}_{i=t+1}^{e-1}(s_{i+1}-c_i)p^i,
        \end{align*}
        which is of the form ${\sum}_{i=0}^rp^{d'_i}$ with $(d'_1,\ldots,d'_r)\le\theta_e$, $d'_r=0$; indeed,
        \begin{align*}
        r&=(s_1-c_0+p)+(p-1)(t-1)-1+{\sum}_{i=t}^{e-1}(s_{i+1}-c_i) \\
        &=p-1-c_0+{\sum}_{i=1}^es_i-{\sum}_{j=t}^{e-1}c_j\le n.
        \end{align*}
        Hence $x'-h-1-p^e\delta_{n-r}\in\Upsilon_x$.
  \item If $s_1=c_0=p-2$ and $c_t=\cdots=c_{e-1}=0$, so that $x'=p^t-2$, then
        \begin{align*}
        x'-h-1=(p^e-p^t)-\Big(&1+p^e+p+{\sum}_{i=1}^{t-1}(p-1)p^i    \\
        &+(s_{t+1}-2)p^t+{\sum}_{i=t+1}^{e-1}s_{i+1}p^i\Big)\in\Upsilon_x,
        \end{align*}
        with the understanding that in $\Lambda_e$ we take $a_t=1$, $a_i=0$ for $i\ne t$, and $$r=p+(p-1)(t-1)+s_{t+1}-2+{\sum}_{i=t+1}^{e-1}s_{i+1}=n.$$
\end{itemize}

Finally, suppose $s_1\le c_0$ and $c_i=s_{i+1}=p-1$ for all $i>0$. In this case, $h=p^e-p+s_1$, $x'=p^e-p+c_0$, and $x=s_1+1-c_0\in\{0,1\}$.
We have
$$(p-2-c_0)p^e-x=(p-c_0)(p^e-1)-\Big(1+p^e+{\sum}_{v=1}^{e}s_vp^{v-1}\Big)\in\Upsilon_x.$$

Therefore, $|\lambda(x)+h+1|\le p^e$ for all $0\le x<p^e$.
By Theorem \ref{thm:main} (b) and Corollary \ref{cor:no-gap},
$$\sigma_0(\tau)=\mu_0(\tau)=\frac{1}{2}\Big(-1-p^e+{\sum}_{v=1}^es_v(p^e-p^{v-1})\Big),$$
recovering \cite{MS19} Theorem 6.2 (a) in the case $p>2$.
\end{exmp}

\ \\
Department of Mathematics, Beijing Technology and Business University, \\
100048, 11\# Fucheng Road, Haidian District, Beijing, China. \\
Haimiao Chen (orcid: 0000-0001-8194-1264)\ \ \ \ chenhm$@$math.pku.edu.cn \\
Yang Li \ \ \ \ lyanghit$@$163.com \\

\end{document}